\documentclass{amsart}
\usepackage{amssymb,euscript,amsmath, mathrsfs, amscd}
\usepackage[dvips]{graphicx}
\usepackage[dvips]{color}

\newcounter{ENUM}

\newcommand{\margh}[1]{}

\input xy
\xyoption{all}
\CompileMatrices

\def\ZZ{{\mathbb Z}}

\def\BB{{\mathbb B}}
\def\PP{{\mathbb P}}
\def\QQ{{\mathbb Q}}

\def\CC{{\mathbb C}}

\def\FF{{\mathbb F}}

\def\cC{{\mathcal C}}

\def\cH{{\mathcal H}}

\def\cX{{\mathcal X}}
\def\cY{{\mathcal Y}}
\def\cZ{{\mathcal Z}}

\def\Aut{\operatorname{Aut}}

\def\Inn{\operatorname{Inn}}

\def\Gal{\operatorname{Gal}}
\def\Ind{\operatorname{Ind}}

\def\Spec{\operatorname{Spec}}

\def\tame{\operatorname{tame}}
\def\PSL{\operatorname{PSL}}
\def\PGL{\operatorname{PGL}}
\def\PGaL{\operatorname{P\Gamma L}}
\newcommand{\adm}{{\rm\scriptstyle adm}}

\newtheorem{thm}{Theorem}[section]
\newtheorem{prop}[thm]{Proposition}
\newtheorem{lem}[thm]{Lemma}
\newtheorem{cor}[thm]{Corollary}

\theoremstyle{definition}
\newtheorem{defn}[thm]{Definition}

\theoremstyle{remark}

\newtheorem{rem}[thm]{Remark}

\numberwithin{equation}{section}
\numberwithin{figure}{section}

\begin{document}
\title{Some $4$-point Hurwitz numbers in positive characteristic}
\author{Irene I. Bouw}
\author{Brian Osserman}
\begin{abstract} 
In this paper, we compute the number of covers of curves with given branch
behavior in characteristic $p$ for one class of examples with
four branch points and degree $p$. Our techniques involve related 
computations in the case of three branch points, and allow us to conclude
in many cases that for a particular choice of degeneration, all the 
covers we consider degenerate to separable (admissible) covers. 
Starting from a
good understanding of the complex case, the proof is centered on the
theory of stable reduction of Galois covers.
\end{abstract}

\maketitle

\section{Introduction}
This paper considers the question of determining the number of covers
between genus-$0$ curves with fixed ramification in positive
characteristic.  More concretely, we consider covers $f:\PP^1\to
\PP^1$ branched at $r$ ordered points $Q_1, \ldots, Q_r$ of fixed {\em
ramification type} $(d; C_1, \ldots, C_r)$, where $d$ is the degree of
$f$ and $C_i=e_1(i)\text{-}\cdots\text{-}e_{s_i}(i)$ is a conjugacy
class in $S_d$. This notation indicates that there are $s_i$
ramification points in the fiber $f^{-1}(Q_i)$, with ramification
indices $e_j(i)$.  The {\em Hurwitz number} $h(d; C_1, \ldots, C_r)$
is the number of covers of fixed ramification type over $\CC$, up to
isomorphism. This number does not depend on the position of the branch
points.  If $p$ is a prime not dividing any of the ramification
indices $e_j(i)$, the {\em $p$-Hurwitz number} $h_p(d; C_1, \ldots,
C_r)$ is the number of covers of fixed ramification type whose branch
points are generic over an algebraically closed field $k$ of
characteristic $p$. The genericity hypothesis is necessary because in
positive characteristic the number of covers often depends on
the position of the branch points.

The only general result on $p$-Hurwitz numbers is that they are always
less than or equal to the Hurwitz number, with equality when the degree
of the Galois closure is prime to $p$. This is because every tame
cover in characteristic $p$ can be lifted  to characteristic $0$,
and in the prime-to-$p$ case, every cover in characteristic $0$ specializes
to a cover in characteristic $p$ with the same ramification type 
(see Corollaire 2.12 of Expos\'e XIII in \cite{sga1}). We say a cover 
has {\em good reduction} when
such a specialization exists. However, in the general case, some covers in 
characteristic $0$ specialize to inseparable covers in characteristic $p$; 
these covers are said to have {\em bad 
reduction}. Thus, the difference $h(d;C_1,\dots,C_r)-h_p(d;C_1,\dots,C_r)$
is the number of covers in characteristic $0$ with generic branch points and 
bad reduction. In \cite{os7} and \cite{os12}, the value 
$h_p(d; e_1, e_2, e_3)$ is computed for genus $0$ covers and any $e_i$ prime 
to $p$ using linear series techniques. In this paper, we treat the 
considerably more difficult case of genus-$0$ covers of type 
$(p; e_1, e_2, e_3, e_4)$. Our main result is the following.

\begin{thm}\label{thm:main} Given $e_1,\dots,e_4$ all less than $p$,
with $\sum_i e_i=2p+2$, we have
\[
h_p(p; e_1, e_2, e_3, e_4)=h(p; e_1, e_2, e_3, e_4)-p.
\]
\end{thm}

An important auxiliary result is the  computation of the
$p$-Hurwitz number $h_p(p; e_1\text{-}e_2, e_3, e_4)$. 

\begin{thm}\label{thm:3-hurwitz} Given odd integers $e_1 , e_2, e_3, e_4 < p$,
with $e_1 + e_2 \leq p$
and $\sum_i e_i=2p+2$,  we have that
\[
h_p(p; e_1\text{-}e_2, e_3, e_4)=
\begin{cases}h(p; e_1\text{-}e_2, e_3, e_4)-(p+1-e_1-e_2): & e_1 \neq e_2, \\
h(p; e_1\text{-}e_2, e_3, e_4)-(p+1-e_1-e_2)/2: & e_1 =e_2. 
\end{cases}
\]
\end{thm}

Corollary \ref{cor:2cyclebad} gives a more general result including
the case that some of the $e_i$ are even, but in some cases we also compute 
the $p$-Hurwitz number only up to a factor $2$. Note that there is an 
explicit formula for $h(p;e_1,e_2,e_3,e_4)$ and 
$h(p;e_1\text{-}e_2,e_3,e_4)$; see Theorem
\ref{hurwitzlem} and Lemma \ref{lem:badtype} below.

Our technique involves the use of ``admissible covers,'' which are
certain covers between degenerate curves (see Section
\ref{sec:char0}).  Admissible covers provide a compactification of the
space of covers of smooth curves in characteristic $0$, but in
positive characteristic this is not the case, and it is an interesting
question when, under a given degeneration of the base, a cover of
smooth curves does in fact have an admissible cover as a limit. In
this case we say the smooth cover has {\em good degeneration}.  In
\cite{bo3} one finds examples of covers with generic branch points
without good degeneration.

In contrast, our technique for proving Theorem \ref{thm:main}
simultaneously shows that many of the examples we consider have good
degeneration. 

\begin{thm}\label{thm:good-degen} Given odd integers
$1< e_1 \leq e_2 \leq e_3 \leq e_4<p$ with $\sum_i e_i=2p+2$, every cover 
of type $(p;e_1,e_2,e_3,e_4)$ with generic branch points $(0,1,\lambda,\infty)$
has good degeneration under the degeneration sending $\lambda$ to $\infty$.
\end{thm}

As with Theorem \ref{thm:3-hurwitz}, our methods do not give a complete
answer in some cases with even $e_i$, but we do prove a more general
result in Theorem \ref{thm:main2}.

Building on the work of Raynaud \cite{ra3}, Wewers uses the theory of
stable reduction in \cite{we1} to give formulas for the number of
covers with three branch points and having Galois closure of
degree strictly divisible by $p$ which have bad reduction to
characteristic $p$. In \cite{b-w6}, some $p$-Hurwitz numbers are
calculated using the existence portion of Wewers' theorems, but these
are in cases which are rigid (meaning the classical Hurwitz number is
$1$) or very close to rigid, so one does not have to carry out
calculations with Wewers' formulas.  In \cite{se6}, Selander uses the
full statement of Wewers' formulas to compute some examples in small
degree. Our result in Theorem  \ref{thm:3-hurwitz} is the first
explicit calculation of an infinite family of $p$-Hurwitz numbers
which fully uses Wewers' formulas, and its proof occupies the bulk of
the present paper.  

We begin in Sections \ref{sec:char0} and \ref{sec:group} by reviewing the 
situation in characteristic $0$ and some group-theoretic background. We
then recall the theory of stable reduction in Section \ref{sec:stable}.
In order to apply Wewers' formulas, in Section \ref{sec:tail} we analyze
the possible structures of the stable reductions which arise, and then
in Section \ref{sec:3pt} we
apply Wewers' formulas to compute the number of smooth covers with a
given stable reduction. Here we are forced to use a trick comparing
the number of covers in the case of interest to the number in a
related case where we know all covers have bad reduction.  
In Section \ref{sec:adm} we then apply Corollary \ref{cor:2cyclebad} as 
well as the formulas
for $h_p(d;e_1,e_2,e_3)$ of \cite{os7} and the classical Hurwitz
number calculations in \cite{o-l2} to estimate the number of
admissible covers in characteristic $p$. This provides a sufficient
lower bound on $h_p(p;e_1,e_2,e_3,e_4)$. Finally, we use the
techniques of \cite{bo4}, again based on stable reduction, to directly
prove in Section \ref{sec:4pt} that $h_p(p;e_1,e_2,e_3,e_4)$ is bounded above 
by $h(p;e_1,e_2,e_3,e_4)-p$. We thus conclude Theorems \ref{thm:main} and
\ref{thm:good-degen}.

We would like to thank Peter M\"uller, Bj\"orn Selander and Robert
Guralnick for helpful discussions.

\section{The characteristic-$0$ situation}\label{sec:char0}

In this paper, we consider covers $f:\PP^1\to \PP^1$ branched at $r$
ordered points $Q_1, \ldots, Q_r$ of fixed {\em ramification type}
$(d; C_1, \ldots, C_r)$, where $d$ is the degree of $f$ and
$C_i=e_1(i)\text{-}\cdots\text{-}e_{s_i}(i)$ is a conjugacy class in
$S_d$. This means that there are $s_i$ ramification points in the
fiber $f^{-1}(Q_i)$, with ramification indices $e_j(i)$.  The {\em
Hurwitz number} $h(d; C_1, \ldots, C_r)$ is the number of covers of
fixed ramification type over $\CC$, up to isomorphism. This number
does not depend on the position of the branch points.

Riemann's Existence Theorem implies that the Hurwitz number $h(d; C_1,
\ldots, C_r)$ is the cardinality of the set of {\em Hurwitz
factorizations} defined as
\[
\{(g_1, \cdots, g_r)\in C_1\times \cdots \times C_r\mid \langle
g_i\rangle\subset S_d\, {\rm transitive },\,  \prod_i g_i=1\}/\sim,
\]
where $\sim$ denotes uniform conjugacy by $S_d$.

The group $\langle g_i \rangle$ is called the {\em monodromy group}
of the corresponding cover. For a fixed monodromy group $G$, a variant 
equivalence relation is given by {\em $G$-Galois covers}, where we work 
with Galois covers together with a fixed isomorphism of the Galois group 
to $G$. 
The group-theoretic interpretation is then that the $g_i$ are in $G$ (with 
the action on a fiber recovered by considering $G$ as a subgroup of 
$S_{|G|}$), and the equivalence relation $\sim_G$ is uniform conjugacy by $G$. 
To contrast with the $G$-Galois case, we sometimes emphasize
that we are working up to $S_d$-conjugacy by referring to the 
corresponding covers as {\em mere covers}.

In this paper, we are mainly interested in the {\em pure-cycle} case,
where every $C_i$ is the
conjugacy class in $S_d$ of a single cycle. In this case, we write
$C_i= e_i$, where $e_i$ is the length of the cycle. A cover $f:Y\to
\PP^1$ over $\CC$ of ramification type $(d; e_1, e_2, \cdots, e_r)$
 has genus $g(Y)=0$ if and only if  $\sum_{i=1}^r e_i=2d-2+r$.

Giving closed formulae for Hurwitz numbers may get very complicated,
even in characteristic zero.  The following result from \cite{o-l2}
illustrates that the genus-$0$ pure-cycle case is more tractable than
the general case, as one may give closed formulae for the Hurwitz numbers, 
at least if the number $r$ of branch points is at most $4$.

\begin{thm}\label{hurwitzlem} Under the hypothesis $\sum_{i=1}^r e_i=2d-2+r$,
we have the following.
\begin{itemize}
\item[(a)] $h(d; e_1, e_2, e_3)=1$.
\item[(b)] $h(d; e_1, e_2, e_3, e_4)=\min_i(e_i(d+1-e_i)).$
\item[(c)] Let $f:\PP^1_\CC\to \PP^1_\CC$ be a cover of ramification
type $(d; e_1, e_2, \ldots, e_r)$ with $r\geq 3$. The Galois group of
the Galois closure of $f$ is either $S_d$ or $A_d$ unless $(d; e_1,
e_2, \ldots, e_r)=(6; 4,4,5)$ in which case the Galois group is $S_5$
acting transitively on $6$ letters.
\end{itemize}
\end{thm}

These statements are Lemma 2.1, Theorem 4.2, and Theorem 5.3 of 
\cite{o-l2}.
We mention that Boccara (\cite{bo5}) proves a partial generalization of 
Theorem \ref{hurwitzlem}.(a). He gives a necessary and sufficient condition
for $h(d; C_1, C_2, \ell)$ to be nonzero in the case that $C_1, C_2$
are arbitrary conjugacy classes of $S_d$ and only $C_3=\ell$ is
assumed to be the conjugacy class of a single cycle.

Later in our analysis we will be required to study covers of
type $(d;e_1\text{-}e_2, e_3,e_4)$, so we mention a result which 
is not stated explicitly in \cite{o-l2},
but which follows easily from the arguments therein. We will only use
the case that $e_4=d$, but we state the result in general since the
argument is the same.

\begin{lem}\label{lem:badtype}
Given $e_1,e_2,e_3,e_4$ and $d$ with $2d+2=\sum_i e_i$ and 
$e_1 + e_2 \leq d$, if $e_1 \neq e_2$ we have
\[
h(d;e_1\text{-}e_2, e_3, e_4)= (d+1-e_1-e_2)\min(e_1,e_2,d+1-e_3,d+1-e_4), 
\]
and if $e_1=e_2$ we have
\[ 
h(d;e_1\text{-}e_2, e_3, e_4)= 
\lceil\frac{1}{2}(d+1-e_1-e_2)\min(d+1-e_3,d+1-e_4)\rceil.
\]

Note that this number is always positive.
In particular, when $e_4=d$ we have
\[
h(d;e_1\text{-}e_2, e_3, d)=
\begin{cases} d+1-e_1-e_2& \text{ if }e_1\neq e_2,\\
(d+2-e_1-e_2)/2&\text{ if }e_1=e_2, d \text{ even},\\
(d+1-e_1-e_2)/2&\text{ if }e_1=e_2, d \text{ odd}.
\end{cases}
\]
\end{lem}

\begin{proof} Without loss of generality, we may assume that
$e_1\leq e_2$ and $e_3 \leq e_4$. Thus, we want to prove that 
$h(d;e_1\text{-}e_2, e_3, e_4)$ is given by the smaller of 
$e_1(d+1-e_1-e_2)$ and $(d+1-e_4)(d+1-e_1-e_2)$
when $e_1\neq e_2$, by $((d+1-e_4)(d+1-e_1-e_2)+1)/2$ when $e_1=e_2$
and all of $d,e_3,e_4$ are even, and by $(d+1-e_4)(d+1-e_1-e_2)/2$ 
otherwise. Even though we do not assume $e_2 \leq e_3$, this formula 
still follows from
the argument of Theorem 4.2.(ii) of \cite{o-l2}.
The first observations to make are that since $e_1+e_2 \leq d$,
we have $e_3+e_4 \geq d+2$, and it follows that although we may not have
$e_2 \leq e_3$, we have $e_1 < e_4$. Moreover, we have $e_1+e_3 \leq d+1$
and $e_2+e_4 \geq d+1$. We are then able to check that the Hurwitz
factorizations $(\sigma_1, \sigma_2, \sigma_3, \sigma_4)$ described in 
case (ii) of {\it loc.\ cit.}\ still give valid Hurwitz factorizations 
$(g_1,g_2,g_3)$ by setting $g_1=\sigma_1\sigma_2$, just as in Proposition
4.7 of {\it loc.\ cit.} Moreover, just as in Proposition 4.9 of 
{\it loc.\ cit.} we find that every Hurwitz factorization must be one
of the enumerated ones. 

It remains to consider when two of the described possibilities yield 
equivalent Hurwitz factorizations. If $e_1\neq e_2$, we can extract
$\sigma_1$ and $\sigma_2$ as the disjoint cycles (of distinct orders)
in $g_1$, so we easily see that the proof of Proposition 4.8 of 
{\it loc. cit.} is still valid. Thus the Hurwitz number is simply the 
number of possibilities enumerated in Theorem 4.2 (ii) of \cite{o-l2}, 
which is the minimum of $e_1(d+1-e_1-e_2)$ and $(d+1-e_4)(d+1-e_1-e_2)$, 
as desired. 

Now suppose $e_1=e_2$. We then check easily that $e_1+e_4 \geq d+1$, so 
that the number of enumerated possibilities is $(d+1-e_4)(d+1-e_1-e_2)$.
Here, we see that we potentially have a given Hurwitz factorization
$(g_1,g_2,g_3)$ being simultaneously equivalent to two of the enumerated
possibilities, since $\sigma_1$ and $\sigma_2$ can be switched. Indeed,
the argument of Proposition 4.8 of {\it loc.\ cit.}\ describing how to 
intrinsically recover the parameters $m,k$ of Theorem 4.2 (ii) of 
{\it loc.\ cit.}\ lets us compute how $m,k$
change under switching $\sigma_1$ and $\sigma_2$, and we find that
the pair $(m,k)$ is sent to $(e_3+2e_4-d-m, e_3+e_4-d-k)$. We thus find
that each Hurwitz factorization is equivalent to two distinct enumerated
possibilities, with the exception that if $d$ and $e_4$ (and therefore
necessarily $e_3$) are even, the Hurwitz factorization corresponding to 
$(m,k)=((e_3+2e_4-d)/2,(e_3+e_4-d)/2)$ is not equivalent to any other.
We therefore conclude the desired statement. 
\end{proof}

We now explain how Theorem 4.2 of \cite{o-l2} can be understood in
terms of degenerations.  Harris and Mumford \cite{h-m2} developed the
theory of {\em admissible covers}, giving a description of the
behavior of branched covers under degeneration.
Admissible covers in the case we are interested in may be described
geometrically as follows: let $X_0$ be the reducible curve consisting
of two smooth rational components $X_0^1$ and $X_0^2$ joined at a
single node $Q$. We suppose we have points $Q_1,Q_2$ on $X_0^1$, and
$Q_3,Q_4$ on $X_0^2$. An {\em admissible cover} of type
$(d;C_1,C_2,\ast,C_3,C_4)$ is then a connected, finite flat cover
$f_0:Y_0 \to X_0$ which is \'etale away from the preimage of $Q$ and
the $Q_i$, and if we denote by $Y_0^1\to X_0^1$ and $Y_0^2\to X_0^2$ the 
(possibly disconnected) covers of $X_0^1$ and $X_0^2$ induced by $f_0$, 
we require also that $Y_0^1\to X_0^1$
has ramification type $(d;C_1,C_2,C)$ for $Q_1,Q_2,Q$ and $Y_0^2\to
X_0^2$ has ramification type $(d;C,C_3,C_4)$ for $Q,Q_3,Q_4$, for some
conjugacy class $C$ in $S_d$, and furthermore that for $P \in f_0^{-1}(Q)$,
the ramification index of $f_0$ at $P$ is the same on $Y_0^1$ and $Y_0^2$.
In characteristic $p$, we further have to require that ramification above 
the node is tame. We refer to $Y_0^1 \to X_0^1$ and $Y_0^2 \to X_0^2$ as 
the {\em component covers} determining $f_0$. When we wish to specify the
class $C$, we say the admissible cover is of type
$(d,C_1,C_2,\ast_C,C_3,C_4)$.

The two
basic theorems on admissible covers concern degeneration and smoothing.
First, in characteristic $0$, or when the monodromy group has order prime to 
$p$, if a family of smooth covers of type $(d;C_1,C_2,C_3,C_4)$ with branch
points $(Q_1,Q_2,Q_3,Q_4)$ is degenerated by sending $Q_3$ to $Q_4$, the
limit is an admissible cover of type $(d;C_1,C_2,\ast,C_3,C_4)$. On the
other hand, given an admissible cover of type $(d;C_1,C_2,\ast,C_3,C_4)$,
irrespective of characteristic there is a deformation to a cover of smooth 
curves, 
which then has type $(d;C_1,C_2,C_3,C_4)$. Such a deformation is not unique
in general; we call the number of smooth covers arising as smoothings of
a given admissible cover (for a fixed smoothing of the base) the 
{\em multiplicity} of the admissible cover.

Suppose we have a family of covers $f:X \to Y$, with smooth generic fiber 
$f_1:X_1 \to Y_1$, and admissible special fiber $f_0:X_0 \to Y_0$.
If we choose local monodromy generators for 
$\pi_1^{\tame}(Y_1 \smallsetminus \{Q_1,Q_2,Q_3,Q_4\})$ which are 
compatible with the degeneration to $Y_0$, we then find that if we have
a branched cover of $Y_1$ corresponding to a Hurwitz factorization
$(g_1,g_2,g_3,g_4)$, the induced admissible cover
of $Y_0$ will have monodromy given by $(g_1,g_2,\rho)$ over
$Y_0^1$ and $(\rho^{-1},g_3,g_4)$ over $Y_0^2$, where
$\rho=g_3 g_4$. The multiplicity of the admissible cover arises
because it may be possible to apply different simultaneous conjugations
to $(g_1,g_2,\rho)$ and to $(\rho^{-1},g_3,g_4)$
while maintaining the relationship between $\rho$ and $\rho^{-1}$.
It is well-known that when $\rho$ is a pure-cycle of order $m$, the
admissible cover has multiplicity $m$, although we recover this fact
independently in our situation as part of the Hurwitz number calculation 
of \cite{o-l2}.

To calculate more generally the multiplicity of an admissible cover of 
the above type, we define the action of the braid operator $Q_3$ on the 
set of Hurwitz factorizations as
\[
Q_3\cdot (g_1, g_2, g_3, g_4)=(g_1g_2g_1g_2^{-1}g_1^{-1}, g_1
g_2g_1^{-1}, g_3, g_4).
\]
One easily checks that $Q_3\cdot\bar{g}$ is again a Hurwitz
factorization of the same ramification type as $\bar{g}$. The multiplicity
of a given admissible cover is the length of the orbit of $Q_3$ acting
on the corresponding Hurwitz factorization.

In this context, we can give the following sharper statement of
Theorem \ref{hurwitzlem} (b), phrased in somewhat different language 
in \cite{o-l2}.

\begin{thm}\label{degenerationlem} 
Given a genus-$0$ ramification type $(d; e_1, e_2, e_3, e_4)$, with
$e_1 \leq e_2 \leq e_3 \leq e_4$
the only possibilities for an admissible cover of type
$(d;e_1,e_2,\ast,e_3,e_4)$ are type 
$(d;e_1,e_2,\ast_{m},e_3,e_4)$ or type
$(d;e_1,e_2,\ast_{e_1\text{-}e_2},e_3,e_4)$.
\begin{itemize}
\item[(a)] Fix $m \geq 1$.  There is at most one admissible cover 
of type $(d;e_1,e_2,\ast_{m},e_3,e_4)$, and if such a cover exists,
it has multiplicity $m$. 
\begin{itemize}
\item[(i)] Suppose that $d+1\leq e_2+e_3$.
There exists an admissible cover of type $(d;e_1,e_2,\ast_{m},e_3,e_4)$
 if and only if
\[
 e_2-e_1+1\leq m\leq 2d+1-e_3-e_4, \qquad m\equiv e_2-e_1+1
\pmod{2}.
\]
\item[(ii)] Suppose that $d+1\geq e_2+e_3$. 
There exists an admissible cover of type $(d;e_1,e_2,\ast_{m},e_3,e_4)$
 if and only if
\[
 e_4-e_3+1\leq m\leq
2d+1-e_3-e_4,\qquad m\equiv e_2-e_1+1 \pmod{2}.
\]
\end{itemize}
\item[(b)] Admissible covers of type 
$(d;e_1,e_2,\ast_{e_1\text{-}e_2},e_3,e_4)$
have multiplicity $1$. The component cover of type 
$(d;e_1,e_2,e_1\text{-}e_2)$ is uniquely determined, so the admissible
cover is determined by its second component cover and the gluing over 
the node. Moreover, the gluing over the node is unique when $e_1 \neq e_2$.
When $e_1=e_2$, there are always two possibilities for gluing except for
a single admissible cover in the case that $e_3,e_4,$ and $d$ are all even.

The number of admissible covers of this type is
\[
\begin{cases}
e_1(d+1-e_1-e_2)& \text{ if }d+1\leq e_2+e_3,\\
(e_3+e_4-d-1)(d+1-e_4)&\text{ if } d+1\geq e_2+e_3.
\end{cases}
\]
\end{itemize}
\end{thm}

\begin{proof} We briefly explain how this follows from Theorem 4.2 of
\cite{o-l2}. As stated above, the possibilities for admissible covers
are determined by pairs $(g_1,g_2,\rho)$, $(\rho^{-1},g_3,g_4)$ where
$(g_1,g_2,g_3,g_4)$ is a Hurwitz factorization of type $(d;e_1,e_2,e_3,e_4)$.
{\it Loc.\ cit.}\ immediately implies that $\rho$ is always either
a single cycle of length $m \geq 1$ or the product of two disjoint
cycles. 

For (a), we find from part (i) of {\it loc.\ cit.}\ that the 
ranges for $m$ (which is $e_3+e_4-2k$ is the notation of {\it loc.\ cit.})
are as asserted, and that for a given $m$,
the number of possibilities with $\rho$ an $m$-cycle is precisely $m$, 
when counted with multiplicity. On the other hand, in this case both
component covers are three-point pure-cycle covers, and thus uniquely
determined (see Theorem \ref{hurwitzlem} (a)). Thus the admissible cover
is unique in this case, with multiplicity $m$. 

For (b), we see by inspection of the description of
part (ii) of {\it loc.\ cit.}\ that $g_1$ is disjoint from $g_2$.
It immediately follows that the braid action is trivial, so the
multiplicity is always $1$, and the asserted count of covers follows
immediately from the proof of Proposition 4.10 of {\it loc.\ cit.} 
Moreover, the component cover of type $(d;e_1,e_2,e_1\text{-}e_2)$
is a disjoint union of covers of type $(e_1;e_1,e_1)$ and $(e_2;e_2,e_2)$
(as well as $d-e_1-e_2$ copies of the trivial cover), so it is uniquely
determined, as asserted. Furthermore, we see that the second component cover
is always a single connected cover of degree $d$, and $g_1,g_2$ are 
recovered as the disjoint cycles of $\rho^{-1}$, so the gluing is unique
when $e_1 \neq e_2$. When $e_1=e_2$, it is possible to swap $g_1$ and
$g_2$, so we see that there are two possibilities for gluing. The 
argument of Lemma \ref{lem:badtype} shows that we do in fact obtain
two distinct admissible covers in this way, except for a single 
cover occurring when $e_3,e_4$ and $d$ are all even.
\end{proof}

\section{Group theory}\label{sec:group}

In several contexts, we will have to calculate monodromy groups other
than those treated by Theorem \ref{hurwitzlem} (c). We will also have
to pass between counting mere covers and counting $G$-Galois covers.
In this section, we give basic group-theoretic results to address
these topics. 

Since we restrict
our attention to covers of prime degree, the following proposition 
will be helpful.

\begin{prop}\label{prop:group}
Let $p$ be a prime number and $G$ a transitive group on $p$
letters. Suppose that $G$ contains a pure cycle of length
$1<e < p-2$. Then $G$ is either $A_p$ or $S_p$. 

Moreover, if $e=p-2$, and $G$ is neither $A_p$ nor $S_p$, then $p=2^r+1$
for some $r$, and $G$ contains a unique
minimal normal subgroup isomorphic to $\PSL_2(2^r)$, and is contained in
$\PGaL_2(2^r)\simeq \PSL_2(2^r)\rtimes \ZZ/r\ZZ$.
If $e=p-1$, and $G$ is not $S_p$, then $G=\FF_p \rtimes \FF_p^*$.
\end{prop}

Note that this does not contradict the exceptional case $d=6$ and
$G=S_5$ in Theorem \ref{hurwitzlem} (c), since we assume that the
degree $d$ is prime.

\begin{proof}
 Since $p$ is prime, $G$ is
necessarily primitive, and a theorem usually attributed to Marggraff 
(\cite{l-t1}) then states that $G$ is at least $(p-e+1)$-transitive. 
When $e\leq p-2$, we have that $p-e+1\geq 3$.
The $2$-transitive permutation groups have been classified by Cameron 
(Section 5 of \cite{ca2}). Specifically, $G$ has a unique minimal 
normal subgroup which is either elementary abelian or one of several
possible simple groups. Since $G$ is at least $3$-transitive, one easily
checks that the elementary abelian case is not possible: indeed, one checks
directly that if a subgroup of a $3$-transitive group inside $S_p$ contains 
an element of prime order $\ell$, then it is not possible for all its
conjugates to commute with one another. 
Similarly, most possibilities
in the simple case cannot be $3$-transitive. If $G$ is not
$S_p$ or $A_p$, then $G$ must have a unique 
minimal normal
subgroup $N$ which is isomorphic to a Mathieu group $M_{11}, M_{23}$,
or to $N=\PSL_2(2^r)$. We then have that $G$ is a subgroup of $\Aut(N)$.
For $N=M_{11},M_{23}$, we have $N=G=\Aut(N)$, and it is easy to check that
the Mathieu groups $M_{11}$ and $M_{23}$ do not contain any single
cycles of order less than $p$, for example with the computer algebra 
package GAP. Therefore
these cases do not occur.  The group $\PSL_2(2^r)$ can only occur if
$p=2^r+1$. In this case, we have that $G$ is a subgroup of
$\Aut(\PSL_2(2^r))= \PGaL_2(2^r)$ and $G$ is at most $3$-transitive, 
so we have $e=p-2$, as desired.

Finally, if $e=p-1$, M\"uller has classified transitive permutation
groups containing $(p-1)$-cycles in Theorem 6.2 of \cite{mu3}, and
we see that the only possibility in prime degree other than $S_p$
is $\FF_p \rtimes \FF_p^*$, as asserted.
\end{proof}

We illustrate the utility of the proposition with:

\begin{cor}\label{3pt-monodromy} Fix $e_1,e_2,e_3,e_4$ with 
$2 \leq e_i \leq p$ for each $i$, and 
$e_1+e_2 \leq p$. For $p>5$, any genus-$0$ cover of type 
$(p;e_1\text{-}e_2,e_3,e_4)$ has monodromy group $S_p$ or $A_p$, with
the latter case occurring precisely when $e_3$ and $e_4$ are odd, and
$e_1+e_2$ is even. For $p=5$, the only exceptional case is type
$(5;2\text{-}2,4,4)$, where the monodromy group is $\FF_5 \rtimes \FF_5^*$.
\end{cor}

\begin{proof} 
Without loss of generality, we assume $e_1\leq e_2$ and $e_3 \leq e_4$.
Applying Proposition \ref{prop:group}, it is clear that
the only possible exception occurs for types with $e_3, e_4 \geq p-2$.
We thus have to treat types $(p;3\text{-}3,p-2,p-2)$, 
$(p;2\text{-}4,p-2,p-2)$, $(p;2\text{-}2,p-2,p)$, $(p;2\text{-}3,p-2,p-1)$, 
and $(p;2\text{-}2,p-1,p-1)$. The fourth case cannot be exceptional since
$G$ contains both a $(p-2)$-cycle and a $(p-1)$-cycle, and the last case
also is ruled out for $p>5$ because $\FF_p \rtimes \FF_p^*$ does not 
contain a $2$-$2$-cycle.

For the first three cases, we must have that $p=2^r+1$ for 
some $r$ and $G$ is a subgroup of $\Gamma:=\PGaL_2(2^r)$. Since $p=2^r+1$ 
is a Fermat prime number, we have that $r$ is a power of $2$. Moreover, 
since $\PSL_2(4)=A_5$ as permutation groups in $S_5$, we may assume 
$r \geq 4$. Since $r$ is even,
any element of order $3$ in $\Gamma \cong \PSL_2(2^r) \rtimes \ZZ/r\ZZ$
must lie inside $\PSL_2(2^r)$, and a non-trivial element of this group
can fix at most $2$ letters. Thus, in order to contain a 
$3\text{-}3$-cycle, we would have to have $6 \leq p=2^r+1 \leq 8$, which 
contradicts the hypothesis $r \geq 4$. This rules out the first case. In 
the second case, if we square the $2$-$4$-cycle we obtain a $2$-$2$-cycle. 
To complete the argument for both the second and third cases it is thus
enough to check directly that if $r>4$, an element of order $2$ cannot
fix precisely $p-4$ letters, ruling out a $2$-$2$-cycle in this case.
It remains only to check directly that $\PGaL_2(16)$ does not contain 
a $2$-$2$-cycle, which one can do directly with GAP.
\end{proof}

Because the theory of stable reduction is developed in the $G$-Galois
context, it is convenient to be able to pass back and forth between
the context of mere covers and of $G$-Galois covers.
The following easy result relates the number of mere covers to the number
of $G$-Galois covers in the case we are interested in.

\begin{lem}\label{Gallem}
Let $f:\PP^1 \to \PP^1$ be a (mere) cover of degree $d$ with monodromy
group $G=A_d$ (respectively, $S_d$). Then the number of $G$-Galois
structures on the Galois closure of $f$ is exactly $2$ (respectively, $1$).
\end{lem}

\begin{proof} 
The case that $G=S_d$ is clear, since conjugacy by $S_d$ is then the
same as conjugacy by $G$. Suppose $G=A_d$, and
let $X=\{(g_1, \ldots, g_r)\mid 
\prod_i g_i=1, \langle g_i\rangle=d\}$. 
Since the centralizer $C_{S_d}(A_d)$ of $A_d$ in $S_d$ is trivial, it
follows that $S_d$ acts freely on $X$, so the number of elements in 
$X_f \subseteq X$ corresponding to $f$ as a mere cover is $|S_d|$. Since 
the center of $G=A_d$ is trivial, $G$ also acts freely on $X$, and $X_f$
breaks into two equivalence classes of $G$-Galois covers, each of
size $|A_d|$.
\end{proof}

\section{Stable reduction}\label{sec:stable} 

In this section, we recall some generalities on stable reduction of
Galois covers of curves, and prove a few simple lemmas as a prelude
to our main calculations. The main references for this section are
\cite{we1} and \cite{bo4}. Since these sources only consider the case
of $G$-Galois covers, we restrict to this situation here as well. Lemma
\ref{Gallem} implies that we may translate results on good or bad
reduction of $G$-Galois covers to results on the stable reduction of the
mere covers, so this is no restriction.

Let $R$ be a discrete valuation ring with fraction field $K$ of
characteristic zero and residue field an algebraically closed field
$k$ of characteristic $p>0$.  Let $f:V=\PP^1_K \to X=\PP^1_K$ be a
degree-$p$ cover branched at $r$ points $Q_1=0,Q_2=1,\ldots,
Q_r=\infty$ over $K$ with ramification type $(p; C_1, \ldots, C_r).$
For the moment, we do not assume that the $C_i$ are the conjugacy
classes of a single cycle.  We denote the Galois closure of $f$ by
$g:Y\to \PP^1$ and let $G$ be its Galois group. Note that $G$ is a
transitive subgroup of $S_p$, and thus has order divisible by $p$.  
Write $H:={\mathop{\rm Gal}}(Y, V)$, a subgroup of index $p$.

We suppose that $Q_i\not\equiv Q_j\pmod{p}$, for $i\neq j$, in other
words, that $(X; \{Q_i\})$ has good reduction as a marked curve. We
assume moreover that $g$ has bad reduction to characteristic $p$, and
denote by $\bar{g}:\bar{Y}\to \bar{X}$ its {\em stable reduction}. The
stable reduction $\bar{g}$ is defined as follows. After replacing $K$
by a finite extension, there exists a unique stable model $\cY$ of
$Y$ as defined in \cite{we1}. We define $\cX=\cY/G$. The stable reduction
$\bar{g}:\bar{Y}:=\cY\otimes_R k\to \bar{X}:=\cX\otimes_R k$ is a
finite map between semistable curves in characteristic $p$; we call
such maps  {\em stable $G$-maps}. We refer to \cite{we1},
Definition 2.14, for a precise definition.

Roughly speaking, the theory of stable reduction proceeds in two steps: 
first, one understands the possibilities for stable $G$-maps, and 
then one counts the number of characteristic-$0$ covers reducing
to each stable $G$-map.

We begin by describing what the stable reduction must look like.
Since $(X; Q_i)$ has good reduction to characteristic $p$, there
exists a model $\cX_0\to \Spec(R)$ such that the $Q_i$ extend to
disjoint sections. There is a unique irreducible component $\bar{X}_0$
of $\bar{X}$, called the {\em original component}, on which the
natural contraction map $\bar{X}\to {\cX}_0\otimes_Rk$ is an
isomorphism. The restriction of $\bar{g}$ to $\bar{X}_0$ is
inseparable. 

 Let $\BB \subseteq \{1,2,\ldots,r\}$ consist of those
indices $i$ such that $C_i$ is not the conjugacy class of a
$p$-cycle. For $i \in \BB$, we have that $Q_i$ specializes to an
irreducible component $\bar{X}_i\neq \bar{X}_0$ of $\bar{X}$.  The 
restriction of $\bar{g}$ to $\bar{X}_i$ is separable, and
$\bar{X}_i$ intersects the rest of $\bar{X}$ in a single point
$\xi_i$. Let $\bar{Y}_i$ be an irreducible component of $\bar{Y}$
above $\bar{X}_i$, and write $\bar{g}_i:\bar{Y}_i\to \bar{X}_i$ for
the restriction of $\bar{g}$ to $\bar{Y}_i$. We denote by $G_i$ the
decomposition group of $\bar{Y}_i$.  We call the components
$\bar{X}_i$ (resp.\ the covers $\bar{g}_i$) for $i\in\BB$ the {\em
  primitive tails} (resp.\ the {\em primitive tail covers}) associated
with the stable reduction. The following definition gives a
characterization of those covers that can arise as primitive tail
covers (compare to \cite{we1}, Section 2.2).

\begin{defn}\label{def:tail} Let $k$ be an algebraically closed field of
 characteristic $p>0$. Let $C$ be a conjugacy class of $S_p$ which
 is not the class  of a $p$-cycle. A {\em primitive tail cover} of
 ramification type $C$ is a $G$-Galois cover $\varphi_C:T_C\to \PP^1_k$
 defined over $k$ which is branched at exactly two points $0, \infty$,
 satisfying the following conditions.
\begin{itemize}
\item[(a)] The Galois group $G_C$ of $\varphi_C$ is a subgroup of
  $S_p$ and contains a subgroup $H$ of index $p$ such that
  $\bar{T}_C:=T_C/H$ has genus $0$.
\item[(b)] The induced map $\bar{\varphi}_C:\bar{T}_C\to \PP^1$ is
  tamely branched at $x=0$, with conjugacy class $C$, and wildly
  branched at $x=\infty$.
\end{itemize}
 If $\varphi$ is a tail cover, we let $h=h(\varphi)$ be
the conductor  and $pm=pm(\varphi)$ the ramification index of
a wild ramification point.

We say that two primitive tail covers $\varphi_i:T_i\to \PP^1_k$ are
{\em isomorphic} if there exists a $G$-equivariant isomorphism $\iota:T_1\to
T_2$. Note that we do not require an isomorphism to send $\bar{T}_1$ to
$\bar{T}_2$.
\end{defn}

Note that an isomorphism $\iota$ of primitive tail covers may be
completed into a commuting square
\[
\begin{CD}
{T}_1 @>{\iota}>> {T}_2\\ 
@V{{\varphi}_1}VV
@VV{{\varphi}_2}V\\ 
\PP^1 @>>>\PP^1.
\end{CD}
\]
Note also that the number of primitive tail covers of fixed
ramification type is finite.

 Since $p$ strictly divides the order of the Galois group $G_C$, we
 conclude that $m$ is prime to $p$. The invariants $h, m$ describe the
 wild ramification of the tail cover $\varphi_C$. The integers $h$ and
 $m$ only depend on the conjugacy class $C$.  In Section
 \ref{sec:tail}, we will show this if $C$ is the class of a single
 cycle or the product of $2$ disjoint cycles, but this holds more
 generally. In the more general set-up of \cite{we1}, Definition 2.9
 it is required that $\sigma:=h/m<1$ as part of the definition of
 primitive tail cover. We will see that in our situation this follows from
 (a).  Moreover, we will show that $\gcd(h, m)=1$ (Lemma
 \ref{lem:tail1}). Summarizing, we find that $(h,m)$ satisfy:
\begin{equation}\label{eq:hm}
m\mid (p-1), \qquad 1\leq h < m, \qquad \gcd(h, m)=1.
\end{equation}

 In the more general set-up of
\cite{we1} there also exists so-called new tails, which satisfy $\sigma>1$. The
following lemma implies that these do not occur in our situation.

\begin{lem}\label{lem:stablered}
The curve $\bar{X}$ consists of at most $r+1$ irreducible components:
the original component $\bar{X}_0$ and primitive tails $\bar{X}_i$ for
all $i \in \BB$. 
\end{lem}

\begin{proof}
In the case that $r=3$ this is proved in \cite{we1}, Section 4.4,
using that the cover is the Galois closure of a genus-$0$ cover of
degree $p$. The general case is a straightforward generalization.
\end{proof}

It remains to discuss the restriction of $\bar{g}$ to the original
component $\bar{X}_0$. As mentioned above, this restriction is
inseparable, and it is described by a so-called deformation datum
(\cite{we1}, Section 1.3).

In order to describe  deformation data, we set some notation. Let 
$\bar{Q}_i$ be the limit on $\bar{X}_0$ of the $Q_i$ for 
$i \not \in \BB$, and set $\bar{Q}_i=\xi_i$ for $i \in \BB$.

\begin{defn}\label{def:dd}
Let $k$ be an algebraically closed field of characteristic $p$.  A
{\em deformation datum} is a pair $(Z, \omega)$, where $Z$ is a smooth
projective curve over $k$ together with a finite Galois cover $g:Z\to
X=\PP^1_k$, and $\omega$ is a meromorphic differential form on $Z$
such that the following conditions hold.
\begin{itemize}
\item[(a)] 
Let $H$ be the Galois group of $Z\to X$.
Then
\[\beta^\ast \omega=\chi(\beta)\cdot
  \omega, \qquad \mbox{for all } \beta\in H.
\]
 Here $\chi:H\to
  \FF_p^\times$ in an injective character.
\item[(b)] The differential form $\omega$ is either logarithmic,
  i.e.\ of the form $\omega={\rm d} f/f$, or exact, i.e.\ of the form
  ${\rm d} f$, for some meromorphic function $f$ on $Z$.
\end{itemize}
\end{defn}

Note that the cover $Z \to X$ is necessarily cyclic.

To a $G$-Galois cover $g:Y\to \PP^1$ with bad reduction, we may
associate a deformation datum, as follows. Choose an irreducible
component $\bar{Y}_0$ of $\bar{Y}$ above the original component
$\bar{X}_0$. Since the restriction $\bar{g}_0:\bar{Y}_0\to\bar{X}_0$
is inseparable and $G\subset S_p$, it follows that
the inertia group $I$ of $\bar{Y}_0$ is cyclic of order $p$, i.e.\ a
Sylow $p$-subgroup of $G$. Since the inertia group is normal in the
decomposition group, the decomposition group $G_0$ of $\bar{Y}_0$ is a
subgroup of $N_{S_p}(I)\simeq \ZZ/p\ZZ\rtimes_\chi\ZZ/p\ZZ^\ast$,
where $\chi:\ZZ/p\ZZ^\ast\to \ZZ/p\ZZ^\ast$ is an injective character.
It follows that the map $\bar{g}_0$ factors as $\bar{g}_0:\bar{Y}_0\to
\bar{Z}_0\to \bar{X}_0$, where $\bar{Y}_0\to \bar{Z}_0$ is inseparable
of degree $p$ and $\bar{Z}_0\to\bar{X}_0$ is separable.  We conclude
that the Galois group $H_0:=\Gal(\bar{Z}_0, \bar{X}_0)$ is a subgroup of
$\ZZ/p\ZZ^\ast\simeq \ZZ/(p-1)\ZZ$. In particular, it follows that
\begin{equation}\label{eq:G0}
G_0\simeq I\rtimes_\chi H_0.
\end{equation}
 The inseparable map $\bar{Y}_0\to
\bar{Z}_0$ is characterized by a differential form $\omega$ on
$\bar{Z}_0$ satisfying the properties of Definition \ref{def:dd}, see
\cite{we1}, Section 1.3.2.  

The differential form $\omega$ is logarithmic if $\bar{Y}_0\to
\bar{Z}_0$ is a ${\boldsymbol \mu}_p$-torsor and exact if this map is
an ${\boldsymbol \alpha}_p$-torsor. A differential form is logarithmic
if and only if it is fixed by the Cartier operator $\cC$ and exact if
and only if it is killed by $\cC$. (See for example \cite{g-s2},
exercise 9.6, for the definition of the Cartier operator and an
outline of these properties.)  Wewers (\cite{we3}, Lemma 2.8)
shows that in the case of covers branched at $r=3$ points the
differential form is always logarithmic.

The deformation datum $(Z, \omega)$ associated to $g$ satisfies the
following compatibilities with the tail covers. We refer to
\cite{we1}, Proposition 1.8 and (2) for proofs of these
statements. For $i\in \BB$, we let $h_i$ (resp.\ $pm_i$) be the
conductor (resp. ramification index) of a wild ramification point of
the corresponding tail cover of type $C_i$, as defined above. We put
$\sigma_i=h_i/m_i$. 
We also use the convention $\sigma_i=0$ for $i \not \in \BB$.
\begin{itemize}
\item[(a)]
 If $C_i$ is the
conjugacy class of a $p$-cycle then $\bar{Q}_i$ is unbranched in
$\bar{Z}_0\to \bar{X}_0$ and $\omega$ has a simple pole at all points
of $\bar{Z}_0$ above $\bar{Q}_i$.
\item[(b)] Otherwise, $\bar{Z}_0\to \bar{X}_0$ is branched of order
  $m_i$ at $\bar{Q}_i$, and $\omega$ has a zero of order $h_i-1$ at
  the points of $\bar{Z}_0$ above $\bar{Q}_i$.
\item[(c)] The map $\bar{Z}_0\to \bar{X}_0$ is unbranched outside
  $\{\bar{Q}_i\}$. All poles and zeros of $\omega$ are above the
  $\bar{Q}_i$.
\item[(d)] The invariants $\sigma_i$ satisfy $\sum_{i\in \BB} \sigma_i=r-2$.
\end{itemize}
The set $(\sigma_i)$ is called the {\em signature} of the deformation
datum $(Z, \omega)$.

\begin{prop}\label{prop:dd} 
Suppose that $r=3, 4$. We fix rational numbers $(\sigma_1, \ldots,
\sigma_r)$ with $\sigma_i\in \frac{1}{p-1}\ZZ$ and $0\leq\sigma_i< 1$,
and $\sum_{i=1}^r \sigma_i=r-2$.  We fix points
$\bar{Q}_1=0, \bar{Q}_2=1, \ldots, \bar{Q}_r=\infty$ on
$\bar{X}_0\simeq \PP^1_k$. Then there exists a deformation datum of
signature $(\sigma_i)$, unique up to scaling. If further the $\bar{Q}_i$
are general, the deformation datum is logarithmic and unique up to 
isomorphism.
\end{prop}

\begin{proof} In the case that $r=3$ this is proved in \cite{we1}.
 (The proof in this case is similar to the proof in the case that
  $r=4$ which we give below.) Suppose that $r=4$. Let $\BB=\{1\leq
  i\leq r\mid \sigma_i\neq 0\}$. We write $\bar{Q}_3=\lambda\in
  \PP^1_k\setminus\{0,1,\infty\}$ and $\sigma_i=a_i/(p-1)$. (If
  $\omega$ is the deformation datum associated with $\bar{g}$, then
  $a_i=h_i (p-1)/m_i.$)

It is shown in \cite{bo4}, Chapter 3, that a deformation datum of
signature $(\sigma_i)$ consists of  a differential form $\omega$ on the cover
$\bar{Z}_0$ of $\bar{X}_0$ defined as a connected component of the
(normalization of the) projective curve with Kummer equation
\begin{equation}\label{eq:Kummer}
z^{p-1}=x^{a_1}(x-1)^{a_2}(x-\lambda)^{a_3}.
\end{equation}
The degree of $\bar{Z}_0\to \bar{X}_0$ is
\[
m:=\frac{p-1}{\gcd(p-1, a_1, a_2, a_3, a_4)}.
\]

The differential form $\omega$ may be written as
\begin{equation}\label{eq:omega}
\omega=\epsilon\frac{z\,{\rm d}
  x}{x(x-1)(x-\lambda)}=
\epsilon \frac{x^{p-a_1}(x-1)^{p-1-a_2}(x-\lambda)^{p-1-a_3}z^p}
{x^p(x-1)^p(x-\lambda)^p}\frac{{\rm d}x}{x},
\end{equation}
where $\epsilon\in k^\times$ is a unit.

To show the existence of the deformation datum, it suffices to show
that one may choose $\epsilon$ such that $\omega$ is logarithmic or
exact, or, equivalently, such that $\omega$ is fixed or killed by the
Cartier operator $\cC$.  It follows from standard properties of the
Cartier operator, (\ref{eq:omega}), and the assumption that
$a_1+a_2+a_3+a_4=2(p-1)$ that $\cC\omega
=c^{1/p}\epsilon^{(1-p)/p}\omega$, where 
\begin{equation}\label{eq:Hasseinv}
c=\sum_{j=\max(0, p-1-a_2-a_4)}^{\min(a_4, p-1-a_3)} {p-1-a_2\choose
      a_4-j}{p-1-a_3\choose j}\lambda^j.
\end{equation}
Note that $c$ is the coefficient of $x^p$ in
$x^{p-a_1}(x-1)^{p-1-a_2}(x-\lambda)^{p-1-a_3}$.  One easily checks
that $c$ is nonzero as  polynomial in $\lambda$. It follows that
$\omega$ defines an exact differential form if and only if $\lambda$
is a zero of the polynomial $c$. This does not happen if
$\{0,1,\lambda, \infty\}$ is general.

We assume that $c(\lambda)\neq 0$.  Since $k$ is algebraically closed,
we may choose $\epsilon\in k^\times$ such that
$\epsilon^{p-1}=c$. Then $\cC\omega=\omega$, and $\omega$ defines a
logarithmic deformation datum. It is easy to see that $\omega$ is
unique, up to multiplication by an element of $\FF_p^\times$.
\end{proof}

\section{The tail covers}\label{sec:tail}
In Section \ref{sec:stable}, we have seen that associated with  a Galois
cover with bad reduction is a set of (primitive) tail covers. In this
section, we analyze the possible tail covers for conjugacy classes
$e\neq p$ and $e_1$-$e_2$ of $S_p$. Recall from Section \ref{sec:char0} that
these are conjugacy classes which occur in the $3$-point covers obtained 
as degeneration of the pure-cycle $4$-point covers.

The following lemma shows the existence of primitive tail covers for the
conjugacy classes occurring in the degeneration of single-cycle
$4$-point covers (Theorem \ref{degenerationlem}).

\begin{lem}\label{lem:tail1}
\begin{itemize}
\item[(a)] Let $2\leq e< p-1$ be an integer. There exists a primitive
  tail cover $\varphi_e:T_e\to \PP^1_k$ of ramification type $e$. Its
  Galois group is $A_p$ if $e$ is odd and $S_p$ if $e$ is even. 
  The wild branch point
  of $\varphi_e$ has inertia group of order $p(p-1)/\gcd(p-1,e-1)
  =:pm_{e}$. The conductor is $h_e:=(p-e)/\gcd(p-1, e-1).$ 
\item[(b)] In the case that $e=p-1$, there exists a primitive tail cover
  $\varphi_e$ of ramification type $e$, 
  with Galois group $\FF_p\rtimes \FF_p^\ast$. The cover is totally
  branched above the wild branch point and has conductor $h_{p-1}=1$.
\item[(c)]
Let $2\leq e_1\leq e_2\leq p-1$ be integers with
  $e_1+e_2\leq p$.  There is a  primitive tail cover
  $\varphi_{e_1,e_2}:T_{e_1,e_2}\to \PP^1_k$ of ramification type
  $e_1$-$e_2$.  The wild branch point of $\varphi_{e_1,e_2}$ has inertia
  group of order $p(p-1)/\gcd(p-1,e_1+e_2-2)=:pm_{e_1, e_2} $. The conductor is
  $h_{e_1,e_2}:=(p+1-e_1-e_2)/\gcd(p-1, e_1+e_2-2).$
\end{itemize}

In all three cases, the tail cover is unique with the given ramification
when considered as a mere cover.
\end{lem}

\begin{proof}
Let $2\leq e\leq p-1$ be an integer.  We define the primitive
tail cover $\varphi_e$ as the Galois closure of the degree-$p$ cover
$\bar{\varphi}_e:\bar{T}_e:=\PP^1\to \PP^1$ given by 
\begin{equation}\label{eq:tail1}
y^p+y^e=x, \qquad (x, y)\mapsto x.
\end{equation}
 One easily checks that this is the unique degree-$p$ cover
  between projective lines with one wild branch point and the required
  tame ramification.  

The decomposition group $G_e$ of $T_e$ is contained in $S_p$. We note
that the normalizer in $S_p$ of a Sylow $p$-subgroup has trivial
center. Therefore the inertia group $I$ of a wild ramification point
of $\varphi_e$ is contained in $\FF_p\rtimes_\chi \FF_p^\ast$, where
$\chi:\FF_p^\ast\to \FF_p^\ast$ is an injective character.  Therefore
it follows from \cite{b-w5}, Proposition 2.2.(i) that $\gcd(h_e,
m_e)=1$. The statement on the wild ramification follows now directly
from the Riemann--Hurwitz formula (as in \cite{we1}, Lemma
4.10). Suppose that $e$ is odd. Then $m_e=(p-1)/\gcd(p-1, e-1)$
divides $(p-1)/2$. Therefore in this case both the tame and the wild
ramification groups are contained in $A_p$. This implies that the
Galois group $G_e$ of $\varphi_e$ is a subgroup of $A_p$.

To prove (a), we suppose that $e\neq p-1$. We show that the Galois
group $G_e$ of $\varphi_e$ is $A_p$ or $S_p$. Suppose that this is not the
case. Proposition \ref{prop:group} implies that
$e=p-2=2^r-1$. Moreover, $G_e$ is a subgroup of $\PGaL_2(2^r)\simeq
\PSL_2(2^r)\rtimes \ZZ/r\ZZ$.  The normalizer in $\PGaL_2(2^r)$ of a
Sylow $p$-subgroup $I$ is $\ZZ/p\ZZ\rtimes \ZZ/2r\ZZ$.  The
computation of the wild ramification shows that the inertia group
$I(\eta)$ of the wild ramification point $\eta$ is isomorphic to
$\ZZ/p\ZZ\rtimes \ZZ/\frac{p-1}{2}\ZZ$.  Therefore $\PGaL_2(2^r)$
contains a subgroup isomorphic to $I(\eta)$ if and only if
$p=17=2^4+1$, in which case $I(\eta)=N_{\PGaL_2(2^r)}(I)$.  We
conclude that if $G_e\not\simeq S_p, A_p$ then $e=15$ and $p=17$.  However,
in this last case one
may check using Magma that a suitable specialization of
(\ref{eq:tail1}) has Galois group $A_{17}$. As before, we conclude
that $G_e\simeq A_{17}$.

 Now suppose that $e=p-1$. It is easy to see that the Galois closure
 of $\bar{\varphi}_{p-1}$ is in this case the cover
 $\varphi_{p-1}:\PP^1\to \PP^1$ obtained by dividing out $\FF_p\rtimes
 \FF_p^\ast\subset \PGL_2(p) = \Aut(\PP^1)$. This proves (b).

Let $e_1, e_2$ be as in the statement of (c).
As in the proof of (a), we define $\varphi_{e_1,
  e_2}$ as the Galois closure of a non-Galois cover
  $\bar{\varphi}_{e_1, e_2}:\bar{T}_{e_1, e_2}\to \PP^1$ of degree
  $p$. The cover $\bar{\varphi}_{e_1, e_2}$, if it exists, is given by
  an equation
\begin{equation}\label{eq:tail2}
F(y):=y^{e_1}(y-1)^{e_2}\tilde{F}(y)=x,\qquad (x, y)\mapsto x,
\end{equation}
where $\tilde{F}(y)=\sum_{i=0}^{p-e_1-e_2} c_i y^i$ has degree
$p-e_1-e_2$. We may assume that $c_{p-e_1-e_2}=1$. The condition that
$\bar{\varphi}_{e_1, e_2}$ has exactly three ramification points
$y=0,1,\infty$ yields the condition $F'(y)=\gamma
t^{e_1-1}(t-1)^{e_2-1}$. Therefore the coefficients of $\tilde{F}$
satisfy the recursion
\begin{equation}\label{eq:rec}
c_i=c_{i-1}\frac{e_1+e_2+i-1}{e_1+i}, \qquad i=1, \ldots, p-e_1-e_2.
\end{equation}
This implies that the $c_i$ are uniquely determined by
$c_{p-e_1-e_2}=1$. Conversely, it follows that the polynomial $F$
defined by these $c_i$ has the required tame ramification. The statement on the
wild ramification follows from the Riemann--Hurwitz formula, as in the
proof of (a).
\end{proof}

It remains to analyze the number of tail covers, and their automorphism
groups. Due to the nature of our argument, we will only need to carry 
out this analysis for the tails of ramification type $e$.
From Lemma \ref{lem:tail1}, it follows already that the map
$\varphi_C:T_C\to \PP^1$ is unique. However, part of the datum of a 
tail cover is an isomorphism $\alpha:~\Gal(T_C, \PP^1)\stackrel{\sim}{\to}
G_C$. For every $\tau\in \Aut(G_C)$, the tuple $(\varphi, \tau\circ
\alpha)$ also defines a tail cover, which is potentially 
non-equivalent. Modification by $\tau$ leaves the cover unchanged as a
$G_C$-Galois cover if and only if $\tau \in \Inn(G_C)$. However, the weaker 
notion of equivalence for tail covers implies that $\tau$ leaves the 
cover unchanged as a tail cover if and only if $\tau$ can be described
as conjugation by an element of $N_{\Aut(T)}(G_C)$. Thus, the number of 
distinct tail covers corresponding to a given mere cover is the order of 
the cokernel of the map 
$$ N_{\Aut(T_C)}(G_C) \to \Aut(G_C) $$
given by conjugation. Denote by
$\Aut_{G_C}(\varphi_C)$ the kernel of this map, or equivalently
the set of $G_C$-equivariant automorphisms of $T_C$. It follows
finally that the number of tail covers corresponding to $\varphi_C$
is
\begin{equation}\label{eq:number-tails}
\frac{|\Aut(G_C)||\Aut_{G_C}(\varphi_C)|}{|N_{\Aut(T_C)}(G_C)|}.
\end{equation}
Finally, denote by $\Aut_{G_C}^0(\varphi_C)\subset\Aut_{G_C}(\varphi_C)$ the 
subset of automorphisms which fix the chosen ramification point $\eta$. We 
now simultaneously compute
these automorphism groups and show that in the single-cycle case, we have
a unique tail cover.

\begin{lem}\label{lem:tail2} Let $2\leq e\leq p-1$ be an integer.
\begin{itemize}
\item[(a)] The group $\Aut_{G_e}(\varphi_e)$
  (resp.\ $\Aut_{G_e}^0(\varphi_e)$) is cyclic of order $(p-e)/2$
  (resp.\ $h_e$) if $e$ is odd and $p-e$ (resp.\ $h_e$) is $e$ is even. 
\item[(b)] There is a unique primitive tail cover of type $e$.
\end{itemize}
\end{lem}

\begin{proof} First note that the definition of $\Aut_{G_e}(\varphi_e)$ 
implies that any element induces an automorphism of any intermediate 
cover of $\varphi_e$, and in particular induces automorphisms of 
$\bar{T}_e$ and $\PP^1$.
Choose a primitive $(p-e)$th root of unity $\zeta\in
\bar{\FF}_p$. Then $\mu(x, y)=(\zeta^ex, \zeta y)$ is an automorphism
of $\bar{T}_e$. One easily checks that $\mu$ generates the group of
automorphisms of $\bar{T}_e$ which induces automorphisms of $\PP^1$
under $\varphi_e$, and that furthermore $T_e$ is Galois over 
$\PP^1/\langle\mu\rangle$, so in particular every element of $\mu$ lifts
to an automorphism of $T_e$. Taking the quotient by the action of $\mu$, 
we obtain a diagram
\begin{equation}\label{tail2eq}
\begin{CD}
\bar{T}_e @>>> \bar{T}'_e=\bar{T}_e/\langle\mu\rangle\\ 
@V{\bar{\varphi}_e}VV
@VV{\bar{\psi}_e}V\\ 
\PP^1 @>>>\PP^1/\langle\mu\rangle.
\end{CD}
\end{equation}

  Since we know the ramification of the other three maps, one easily
  computes that the tame ramification of $\bar{\psi}_e$ is
  $e$-$(p-e)$. Let $\psi_e: T_e'\to \PP^1$ be the Galois closure of
  $\bar{\psi}_e$. 

We now specialize to the case that $e$ is odd. Since 
  $G_e=A_p$ does not contain an element of cycle type $e$-$(p-e)$, it follows 
  that the Galois group $G'$ of $\psi_e$ is $S_p$. 
Therefore it follows by degree considerations that the cover 
$T_e \to T_e'$ is cyclic of degree $(p-e)/2$.  
Denote by $Q$ the Galois group of the cover $T_e\to \PP^1/\langle\mu\rangle$. 
This is a group of order
$p!(p-e)/2$, which contains normal subgroups isomorphic to $A_p$ and
$\ZZ/\frac{p-e}{2}\ZZ$, respectively.  
It follows that $Q = \ZZ/\frac{p-e}{2}\ZZ \rtimes S_p$.
Note that $\Aut_{G_e}(\varphi_e)$ is necessarily a subgroup of
$Q$. In fact, it is precisely the subgroup of $Q$ which commutes
with every element of $A_p \subseteq Q$. One easily checks that the
semidirect product cannot be split, and that $\Aut_{G_e}(\varphi_e)$
is precisely the normal subgroup $\ZZ/\frac{p-e}{2}\ZZ$, that is the
Galois group of $T_e$ over $T'_e$.

To compute $\Aut^0_{G_e}(\varphi_e)$ we  need to compute the
order of the inertia group  of a wild ramification point of 
$T_e$ in the map $T_e\to T'_e$.
 Since a wild ramification point of $T_e'$
  has inertia group of order $p(p-1)=pm_e\gcd(p-1, e-1)$, we know the
orders of the inertia groups of three of the four maps, and conclude
  that $\Aut^0_{G_e}(\varphi_e)$ has order $h_e=(p-e)/\gcd(p-1, e-1)$.
This proves (a) in the case $e$ is odd. 

For (b), we simply observe that since $Q \subset N_{\Aut(T_e)}(G_e)$,
we have
$$\frac{|\Aut(G_e)||\Aut_{G_e}(\varphi_e)|}{|N_{\Aut(T_e)}(G_e)|} \leq
\frac{p!\frac{p-e}{2}}{|Q|}=1,$$
so the tail cover is unique, as desired.

We now treat the case that $e$ is even. For (a), if $e<p-1$, the Galois 
group of $\bar{\psi}_e$ is equal to the Galois group of $\bar{\varphi}_e$, 
which is isomorphic to $S_p$.  We conclude that the degree of 
$T_e\to T'_e$ is $p-e$ in this case, and the group $Q$ defined as above 
is a direct product $\ZZ/(p-e)\ZZ \times S_p$.
Similarly to the case that $e$ is odd, we conclude that
$\Aut_{G_e}(\varphi_e)$ (resp.\ $\Aut_{G_e}^0(\varphi_e)$) is cyclic
of order $p-e$ (resp.\ $h_e$) in this case, as desired. 
On the other hand, if $e=p-1$, we have 
that $p-e=1$, hence $\mu$ is trivial, and we again conclude that (a) holds.
Finally, (b) is trivial: if $e<p-1$, the Galois group of $\varphi_e$ is 
$S_p$ and $\Aut(S_p)=S_p$. Therefore there is a unique tail cover in this 
case. The same conclusion holds in the case that $e=p-1$, since 
$G_{p-1}\simeq \FF_p\rtimes_\chi\FF_p^\ast$ and $\Aut(G_{p-1})=G_{p-1}$.
The statement of the lemma follows.
\end{proof}

\begin{rem} In the case of $e_1$-$e_2$ tail covers, there may in fact
be more than one structure on a given mere cover. However, we will
not need to know this number for our argument.
\end{rem}

\section{Reduction of $3$-point covers}\label{sec:3pt}
In this section, we (almost) compute the number of $3$-point covers
with bad reduction for ramification types
$(p;e_1\text{-}e_2,e_3,e_4)$. More precisely, we compute this number
in the case that not both $e_3$ and $e_1+e_2$ are even. In the remaining 
case, we only compute this number up to a factor $2$, which is good enough
for our purposes. Although we restrict to types of the above form, 
our strategy applies somewhat more generally.
The results of this section rely on the results of Wewers \cite{we1},
who gives a precise formula for the number of lifts of a given special
$G$-map (Section \ref{sec:stable}) in the $3$-point case. 

We fix a type $\tau=(p;e_1\text{-}e_2,e_3,e_4)$ satisfying the
genus-$0$ condition $\sum_i e_i=2p+2$.  We allow $e_3$ or $e_4$ to be
$p$, although this is not the case that ultimately interests us; see
below for an explanation. We do however assume throughout that we are
not in the exceptional case $\tau=(5;2\text{-}2,4,4)$. According to
Lemma \ref{lem:tail1}, we may fix a set of primitive tail covers
$\bar{g}_i$ of type $C_i$, for $i$ such that $C_i\neq p$. Moreover,
by Proposition \ref{prop:dd} we have a (unique) deformation datum,
so we know there exists at least one special $G$-map $\bar{g}$ of type 
$\tau$. Lemma \ref{lem:tail2} implies moreover that the number of special 
$G$-maps is equal to the number of $e_1$-$e_2$ tail covers. 
Wewers (\cite{we1}, Theorem 3) shows that there exists a $G$-Galois cover
$g:Y\to \PP^1$ in characteristic zero with bad reduction to
characteristic $p$, and more specifically with stable reduction equal to
the given special $G$-map $\bar{g}$.
Moreover, Wewers gives a formula for the number 
$\tilde{L}(\bar{g})$ for lifts of the given special $G$-map $\bar{g}$. 

In order to state his formula, we need to introduce one more invariant.
Let $\Aut_G^0(\bar{g})$ be the group of $G$-equivariant automorphisms of
$\bar{Y}$ which induce the identity on the original component
$\bar{X}_0$. 
Choose $\gamma\in \Aut_G^0(\bar{g})$, and consider the restriction of
$\gamma$ to the original component $\bar{X}_0$. Let $\bar{Y}_0$ be an
irreducible component of $\bar{Y}$ above $\bar{X}_0$ whose inertia
group is the fixed Sylow $p$-subgroup $I$ of $G$. As in (\ref{eq:G0}),
we write $G_0=I\rtimes_\chi H_0\subset \FF_p\rtimes_\chi \FF_p^\ast$
for the decomposition group of $\bar{Y}_0$.  Wewers (\cite{we1}, proof
of Lemma 2.17) shows that $\gamma_0:=\gamma|_{\bar{Y}_0}$ centralizes
$H_0$ and normalizes $I$, i.e.\ $\gamma_0\in C_{N_G(I)}(H_0)$. Since
$\bar{Y}|_{\bar{X}_0}=\Ind_{G_0}^G\bar{Y}$ and $\gamma$ is
$G$-equivariant, it follows that the restriction of $\gamma$ to
$\bar{X}_0$ is uniquely determined by $\gamma_0$. We denote by $n'(\tau)$
the order of the subgroup consisting of those  $\gamma_0\in C_{N_G(I)}(H_0)$
such that there exists a $\gamma\in \Aut_G^0(\bar{g})$ with
$\gamma|_{\bar{Y}_0}=\gamma_0$. Our notation is justified by Corollary 
\ref{cor:nprime-defnd} below.

Wewers (\cite{we1}, Corollary 4.8) shows that
\begin{equation}\label{eq:pd}
|\tilde{L}(\bar{g})|=\frac{p-1}{n'(\tau)}\prod_{i\in
  \BB}\frac{h_{C_i}}{|\Aut_{G_{C_i}}^0(\bar{g}_{C_i})|}.
 \end{equation}
The numbers are as defined in Section \ref{sec:stable}.  (Note that the
group $\Aut_{G_{C_i}}^0(\bar{g}_{C_i})$ is defined differently from
the group $\Aut_G^0(\bar{g})$.)

To compute the number of curves with bad reduction, we need to
compute the number $n'(\tau)$ defined above. As explained by Wewers
(\cite[Lemma 2.17]{we1}), one may express the number $n'(\tau)$ in terms of
certain groups of automorphisms of the tail covers. However, there is a
mistake in the concrete description he gives of $\Aut_G^0(\bar{g})$ in
terms of the tails, therefore we do not use Wewers' description. For a
corrected version, we refer to the manuscript \cite{se6}.

The difficulty we face in using Wewers' formula directly is that
we do not know the Galois group $G_{e_1\text{-}e_2}$ of the $e_1$-$e_2$
tail. This prevents us from directly computing both the number of 
$e_1$-$e_2$ tails, and the term $n'(\tau)$.
We avoid this problem by using the following trick.
 We first consider covers of type $\tau^\ast=(p;
 e_1\text{-}e_2, \varepsilon, p)$, with $\varepsilon=p+2-e_1-e_2$, which 
 all have bad reduction. This
 observation lets us compute $n'(\tau^\ast)$ from Wewers' formula.
 We then show that for covers of type $\tau=(p;
 e_1\text{-}e_2,e_3,e_4)$, the number $n'(\tau)$ essentially only depends on
 $e_1$ and $e_2$, allowing us to compute $n'(\tau)$ from $n'(\tau^\ast)$.
 A problem with this method is that in the case that the Galois
 groups of covers with type $\tau$ and $\tau^\ast$ are not equal, the
 numbers $n'(\tau)$ and $n'(\tau^\ast)$ may 
 differ by a factor $2$. Therefore in this case, we are able to
 determine the number of covers of type $\tau$ with bad reduction only
 up to a factor $2$.

  In Lemma \ref{lem:badtype}, we have counted non-Galois covers, but in 
this section, we deal with Galois covers. 
 Let $G(\tau)$ be the Galois group of a cover of type $\tau$. This
 group is well-defined and either $A_p$ or $S_p$, by Corollary
 \ref{3pt-monodromy}.  We write $\gamma(\tau)$ for the quotient of the
 number of Galois covers of type $\tau$ by the Hurwitz number
 $h(\tau)$. By Lemma \ref{Gallem}, it follows that $\gamma(\tau)$
 is $2$ if $G$ is $A_p$ and $1$ if it is $S_p$.  The number
 $\gamma(\tau)$ will drop out from the formulas as soon as we
 pass back to the non-Galois situation in Section \ref{sec:adm}.

We first compute the number
$n'(\tau^\ast)$. We note that by Corollary \ref{3pt-monodromy}, 
the Galois group $G(\tau^\ast)$ of a
cover of type $\tau^\ast$ is $A_p$ if $e_1+e_2$ is even and $S_p$
otherwise. In particular, we see that $G(\tau)=G(\tau^\ast)$ unless
$e_1+e_2$ and $e_3$ are both even. In this case we have that
$G(\tau)=S_p$ and $G(\tau^\ast)=A_p$. Recall from Lemma \ref{lem:tail2} that
there is a unique
tail cover for the single-cycle tails. We denote by $N_{e_1\text{-}e_2}$
the number of $e_1\text{-}e_2$ tails, and by $\Aut^0_{e_1\text{-}e_2}$
the group $\Aut^0_{G_{e_1\text{-}e_2}}(\bar{g}_{e_1\text{-}e_2})$ for
any tail cover $\bar{g}_{e_1\text{-}e_2}$ as in Lemma \ref{lem:tail1}.
Note that since $\bar{g}_{e_1\text{-}e_2}$ is unique as a mere cover,
and the definition of $\Aut^0_{G_{e_1\text{-}e_2}}(\bar{g}_{e_1\text{-}e_2})$
is independent of the $G$-structure, this notation is well-defined.
We similarly have from \eqref{eq:pd} that $|\tilde{L}(\bar{g})|$
depends only on $\tau$, so we write $\tilde{L}(\tau):=|\tilde{L}(\bar{g})|$
for any special $G$-map $\bar{g}$ of type $\tau$.

\begin{lem}\label{lem:tauast}
Let $\tau^\ast$ be as above. Then
\[
n'(\tau^\ast)= \frac{(1+\delta_{e_1,e_2})N_{e_1\text{-}e_2}(p-1)}
{\gcd(p-1, e_1+e_2-2)\gamma(\tau^\ast)|\Aut^0_{e_1\text{-}e_2}|}.
\]
Here $\delta_{e_1,e_2}$ is the Kronecker $\delta$.
\end{lem}

\begin{proof}
 Lemma \ref{lem:badtype} implies that the Hurwitz number $h(\tau^\ast)$
 equals $(p+1-e_1-e_2)/2$ if $e_1=e_2$ and $(p+1-e_1-e_2)$
 otherwise. Since all covers of type $\tau^\ast$ have bad reduction,
 $h(\tau^\ast)\gamma(\tau^\ast)$ is equal to 
 $N_{e_1\text{-}e_2} \cdot \tilde{L}(\tau^\ast)$. The statement of the 
 lemma follows by applying Lemmas \ref{lem:tail1}.(c), \ref{lem:tail2},
  and (\ref{eq:pd}).
\end{proof}

We now analyze $n'$ in earnest. For convenience, for $i\in \BB$ we also 
introduce the notation $\widetilde{\Aut}_{G_i}(\bar{g}_i)$
for the group of $G$-equivariant automorphisms of the induced tail cover
$\Ind_{G_i}^G(\bar{g}_i)$. Recall also that $\xi_i$ is the node connecting
$\bar{X}_0$ to $\bar{X}_i$. We note that $n'$ may be analyzed tail by tail,
in the sense that given $\gamma_0 \in C_{N_G(I)}(H_0)$, we have that
$\gamma_0$ lifts to $\Aut^0(\bar{g})$ if and only if for each $i \in \BB$,
there is some $\gamma_i \in \widetilde{\Aut}_{G_i}(\bar{g}_i)$ whose
action on $\bar{g}_i^{-1}(\xi_i)$ is compatible with $\gamma_0$.
The basic proposition underlying the behavior of $n'$ is then the following:

\begin{prop}\label{prop:nprime-desc} Suppose $G=S_p$ or $A_p$, and we have a 
special $G$-map $\bar{g}:\bar{Y}\to\bar{X}$. Then:
\begin{itemize}
\item[(a)] For $i\in \BB$, the 
$G$-equivariant automorphisms of $\bar{g}^{-1}(\xi_i)$ form a cyclic group. 

\item[(b)] Given an element $\gamma_0 \in C_{N_G(I)}(H_0)$ and $i \in \BB$,
there exists $\gamma_i\in \widetilde{Aut}_{G_i}(\bar{g}_i)$ 
agreeing with the action of $\gamma_0$ on $\bar{g}^{-1}(\xi_i)$ if and only
if there exists $\gamma_i'\in \widetilde{Aut}_{G_i}(\bar{g}_i)$ 
having the same orbit length on $\bar{g}^{-1}(\xi_i)$ as $\gamma_0$ has.
\end{itemize}
\end{prop}

\begin{proof} For (a), if $\tilde{\xi}_i$ is a point above
$\xi_i$ lying on the chosen component $\bar{Y}_0$, one easily checks 
that a $G$-equivariant automorphism $\gamma$ of $\bar{g}^{-1}(\xi_i)$ is
determined by where it sends $\tilde{\xi}_i$, which can in turn be
represented by an element $g \in G$ chosen so that 
$g(\tilde{\xi}_i)=\gamma(\tilde{\xi}_i)$. Note that $\gamma \neq g$;
in fact, if $\gamma,\gamma'$ are determined by $g,g'$, the composition
law is that $\gamma \circ \gamma'$ corresponds to $g' g$.
Such a $g$ yields a choice of $\gamma$ if and only if we have
the equality of stabilizers $G_{\tilde{\xi}_i}=G_{g(\tilde{\xi}_i)}$.
Now, any $h \in G_{\tilde{\xi}_i}$ is necessarily in $G_0$, 
and using that $I \subseteq G_{\tilde{\xi}_i}$, we find that we must have
$g I g^{-1} \subseteq G_0$. But $I$ contains the only $p$-cycles in
$G_0$, so we conclude $g\in N_G I$. However, since $I$ fixed
$\tilde{\xi}_i$, the choices of $G$ may be taken modulo $I$, so we
conclude that they lie in $N_G I/I$. Finally, since $G=S_p$ or $A_p$,
we have that $N_G I/I$ is cyclic, isomorphic to $\ZZ/(p-1)\ZZ$ if
$G=S_p$ and to $\ZZ/(\frac{p-1}{2})\ZZ$ if $G=A_p$.

(b) then follows immediately, since the actions of both $\gamma_0$ and 
$\gamma_i'$ on $\bar{g}^{-1}(\xi_i)$ lie in the same cyclic group; we
can take $\gamma_i$ to be an appropriate power of $\gamma_i'$. 
\end{proof}

\begin{cor}\label{cor:nprime-defnd} For $\tau$ as above, $n'(\tau)$ is 
well defined.
\end{cor}

\begin{proof} We know that $G=S_p$ or $A_p$, and we also know by
Proposition \ref{prop:dd} and Lemma \ref{lem:tail1} that the deformation 
datum is uniquely determined, and so
are the tail covers, at least as mere covers. But the description of
$n'(\tau)$ given by Proposition \ref{prop:nprime-desc} is visibly
independent of the $G$-structure on the tail covers, so we obtain
the desired statement.
\end{proof}

We can now obtain the desired comparison of $n'(\tau)$ with $n'(\tau^\ast)$.

\begin{prop}\label{prop:n'}
Let $\tau=(p; e_1$-$e_2, e_3, e_4)$ be a type satisfying the genus-$0$
condition, and let $\tau^\ast$ be the corresponding modified
type.
Then if $G(\tau)=G(\tau^\ast)$ we have
  $n'(\tau)=n'(\tau^\ast)$. Otherwise, $n'(\tau)\in \{2n'(\tau^\ast),
  n'(\tau^\ast)\}$.
\end{prop}

\begin{proof} 
Let $\gamma_0$ be a generator of $C_{N_G(I)}(H_0)$. 
We ask which powers of $\gamma_0$ extend to an element of $\Aut_G^0(\bar{g})$,
and we analyze this question tail by tail. Fix a tail $\bar{X}_i$,
and suppose that it is a single-cycle tail of length $e:=e_i$. The crucial 
assertion is that $\gamma_0$ itself (and hence all its powers) always 
extends to $\bar{X}_i$.

First suppose that $e<p-1$ is even. Thus $G=G_i=S_p$, and
$\widetilde{\Aut}_G(\bar{g}_i)=\Aut_{G_i}(\bar{g}_i)$. Now, $\gamma_0$ acts
on the fiber of $\xi_i$ with orbit length $(p-1)/m_e=\gcd(p-1,p-e)$. 
On the other hand, by Lemma \ref{lem:tail1} we have that 
$h_e=(p-e)/\gcd(p-1,e-1)$. Lemma \ref{lem:tail2} implies 
that if $\gamma_i \in \Aut_{G_i}(\bar{g}_i)$ 
is a generator, then the order of $\gamma_i$ is $p-e$, and also that
$\Aut^0_{G_i}(\bar{g}_i)$ has order $h_e$. We conclude that an orbit of
$\gamma_i$ has length $\gcd(p-1,e-1)=\gcd(p-1,p-e)$, and thus by Proposition 
\ref{prop:nprime-desc} that $\gamma_0$ extends to $\bar{X}_i$, as claimed. 

The next case is that $e$ is odd, and $G=A_p$.
This proceeds exactly as before, except that both orbits in question have
length $\gcd(p-1,p-e)/2$.
Now, suppose $e$ is odd, but $G=S_p$. Then the orbit length of
$\gamma_0$ is $\gcd(p-1,p-e)$. 
We have $\widetilde{\Aut}_{G}(\bar{g}_i)$ equal to the $G$-equivariant
automorphisms of $\Ind_{A_p}^{S_p}(\bar{g}_i)$. These contain induced
copies of the $G$-equivariant automorphisms of $\bar{g}_i$, so in
particular we know we have elements of orbit length 
$\gcd(p-1,e-1)/2$. However, in fact one also has a $G$-equivariant
automorphism exchanging the two copies of $\bar{g}_i$, and whose square
is the generator of the $A_p$-equivariant automorphisms of $\bar{g}_i$. One 
may think of this as coming from the automorphism constructing in
Lemma \ref{lem:tail2} inducing the isomorphism 
between the two different $A_p$-structures on the tail cover. We thus
have an element of $\widetilde{\Aut}_{G}(\bar{g}_i)$ of orbit length
$\gcd(p-1,e-1)$, and $\gamma_0$ extends to the tail in this case as well.

Finally, if $e=p-1$ then $m_i=p-1$ and thus $\gamma_0$ acts as the
identity on the fiber of $\xi_i$. The
claim is trivially satisfied in this case. 

It follows that extending $\gamma_0$ to the $e$-tails imposes no
condition when $e<p$, and of course we do not have tails in the case
that $e=p$. Therefore the only non-trivial condition imposed in extending
$\gamma_0$ is the extension to the $e_1$-$e_2$-tail.

In the case that $G(\tau)=G(\tau^\ast)$ we conclude the desired statement
from Proposition \ref{prop:nprime-desc}, since the orbit lengths in
question are clearly the same in both cases. Suppose that $G(\tau)\neq
G(\tau^\ast)$. This happens if and only if both $e_1+e_2$ and $e_3$
are even. In this case we have that $G(\tau)=S_p$ and
$G(\tau^\ast)=A_p$. Here, we necessarily have that $e_1+e_2,e_3,e_4$ are all
even, so the only conditions imposed on either $n'(\tau)$ or $n'(\tau^\ast)$
come from the $e_1$-$e_2$ tail. Since the orbit of $\gamma_0$ is twice as
long in the case of $\tau$, the answers can differ by at most a factor of
$2$ in this case, as desired.
\end{proof}

Let $2\leq e_1\leq e_2\leq e_3\leq e_4<p$ be integers with $\sum_i
e_i=2p+2$ and $e_1+e_2\leq p$.  The following corollary translates
Proposition \ref{prop:n'} into an estimate for the number of Galois
covers of type $\tau=(p; e_1$-$e_2, e_3, e_4)$ with bad
reduction. Theorem \ref{thm:3-hurwitz} is a special case.

\begin{cor}\label{cor:2cyclebad} Let $\tau=(p; e_1\text{-}e_2, e_3, e_4)$ 
with $\tau \neq (5;2\text{-}2,4,4)$. The number of mere covers of type $\tau$
with bad reduction to characteristic $p$ is equal to
\[
\begin{cases}
\delta(\tau)(p+1-e_1-e_2)&\text{ if }e_1\neq e_2,\\
\delta(\tau)(p+1-e_1-e_2)/2&\text{ if } e_1=e_2,
\end{cases}
\]
where $\delta(\tau)\in \{1,2\}$, and $\delta=1$ unless $e_1+e_2$ and $e_3$
are both even.
\end{cor}

\begin{proof}
We recall that the number of Galois covers of type $\tau$ with bad
reduction is equal to $N_{e_1\text{-}e_2} \cdot \tilde{L}(\tau^\ast)$. 
It follows from Lemma
\ref{lem:tail1}.(c) and (\ref{eq:pd}) 
that this number is
\[
\begin{cases}
\displaystyle{
  \frac{\gamma(\tau^\ast)n'(\tau^\ast)}{n'(\tau)}(p+1-e_1-e_2)}&
  \text{ if }e_1\neq e_2,\\
  \displaystyle{\frac{\gamma(\tau^\ast)n'(\tau^\ast)}{n'(\tau)}(p+1-e_1-e_2)/2}&
  \text{ if }e_1= e_2.
\end{cases}
\]
 The definition of the Galois factor $\gamma(\tau)$ implies that the
 number of mere covers of type $\tau$ with bad reduction is
\[
\begin{cases}
\displaystyle{\frac{\gamma(\tau^\ast)}{\gamma(\tau)}
  \frac{n'(\tau^\ast)}{n'(\tau)}(p+1-e_1-e_2)}& \text{ if }e_1\neq
e_2,\\ \displaystyle{\frac{\gamma(\tau^\ast)}{\gamma(\tau)}
  \frac{n'(\tau^\ast)}{n'(\tau)}(p+1-e_1-e_2)/2}& \text{ if }e_1= e_2.
\end{cases}
\]

 Proposition \ref{prop:n'} implies that $n'(\tau)/n'(\tau^\ast)\in
 \{1,2\}$, and is equal to $1$ unless $e_1+e_2,e_3,e_4$ are all even.
 Moreover, if $n'(\tau)\neq n'(\tau^\ast)$ then
 $\gamma(\tau^\ast)/\gamma(\tau)=2$. 
The statement of the corollary follows from this.
\end{proof}

\begin{rem} Similar to the proof of Corollary \ref{cor:2cyclebad}, one may
 show that every genus-$0$ three-point cover of type $(p; e_1, e_2,
 e_3)$ has bad reduction. We do not include this proof here, as a
 proof of this result using linear series already occurs in
 \cite{os7}, Theorem 4.2. 
\end{rem}

\section{Reduction of admissible covers}\label{sec:adm}

In this section, we return to the case of non-Galois covers, and use the 
results of Section \ref{sec:3pt} to compute the number of 
``admissible covers with good reduction''. We start by defining what we
mean by this. As always, we fix a type $(p; e_1, e_2, e_3, e_4)$ with
$1<e_1\leq e_2\leq e_3\leq e_4<p$ satisfying the genus-$0$ condition
$\sum_i e_i=2p+2$. As in Section \ref{sec:char0}, we consider
admissible degenerations of type $(p;e_1, e_2,\ast,e_3, e_4)$, which
means that $Q_3=\lambda\equiv Q_4=\infty\pmod{p}$.  Recall from
Section \ref{sec:char0} that in positive characteristic not every smooth 
cover degenerates to an admissible cover, as a degeneration might
become inseparable.  The number of admissible covers (even counted
with multiplicity) is still bounded by the number of smooth covers,
but equality need not hold.

\begin{defn}\label{def:phadm} 
We define $h^\adm_p(p; e_1, e_2,\ast,e_3, e_4)$ as the number of admissible
covers of type $(p; e_1,e_2, \ast, e_3,e_4)$, counted with multiplicity, over
an algebraically closed field of characteristic $p$.
\end{defn}

The following proposition is the main result of this
section.

\begin{prop}\label{prop:admbad} 
The assumptions on the type $\tau=(p;e_1, e_2, e_3, e_4)$ are as above. Then
\[
h^\adm_p(p; e_1, e_2,\ast, e_3, e_4)>h(p; e_1, e_2, e_3, e_4)-2p,
\]
and
\[
h^\adm_p(p; e_1, e_2,\ast, e_3, e_4)=h(p; e_1, e_2, e_3, e_4)-p
\]
unless $e_1+e_2$ and $e_3$ are both even.
\end{prop}

\begin{proof} 
We begin by noting that in the case $\tau=(5;2,2,4,4)$ corresponding to
the exceptional case of Corollary \ref{3pt-monodromy},
the assertion of the proposition is automatic since $h(5;2,2,4,4)=8<10$.
We may therefore assume that $\tau \neq (5;2,2,4,4)$.

We use the description of the admissible covers in characteristic
 zero (Theorem \ref{degenerationlem}) and the results of Section
 \ref{sec:3pt} to estimate the number of admissible covers with good
 reduction to characteristic $p$, i.e.\  that remain separable.

We first consider the pure-cycle case, i.e.\ the case of Theorem
\ref{degenerationlem}.(a). Let $m$ be an integer satisfying the
conditions of {\it loc.\ cit.} We write ${f}_0:{V}_0\to {X}_0$ for the
corresponding admissible cover. Recall from Section \ref{sec:char0} that
$\bar{X}$ consists of two projective lines ${X}^1_0, {X}^2_0$
intersecting in one point. Choose an irreducible component ${Y}^i_0$
of ${Y}_0$ above ${X}^i_0$, and write ${f}^i_0:{Y}^i_0\to {X}^i_0$ for
the restriction. These are covers of type $(d_1; e_1, e_2, m)$ and
$(d_2; m, e_3, e_4)$ with $d_i\leq p$, respectively.  The
admissible cover ${f}_0$ has good reduction to characteristic $p$ if
and only if both three-point covers ${f}^i_0$  have good
reduction.

It is shown in \cite{os7}, Theorem 4.2,
that a genus-$0$ three-point cover of type $(d;a,b,c)$ with $a,b,c<p$
has good reduction to characteristic $p$ if and only if its degree $d$
is strictly less than $p$. 
 Since the degree $d_2$ of the cover ${f}^2_0$ is always at least as
 large as the other degree $d_1$, it is enough to calculate when
 $d_2<p$. The Riemann--Hurwitz formula implies that
 $d_2=(m+e_3+e_4-1)/2$. Therefore the condition $d_2<p$ is
 equivalent to the inequality
$$e_3+e_4+m \leq 2p-1.$$

Since we assumed the existence of an admissible cover with $\rho$ an 
$m$-cycle, it follows from
Theorem \ref{degenerationlem}.(a) that $m \leq
2d+1-e_3-e_4=2p+1-e_3-e_4$. We find that $d_2<p$ unless
$m=2p+1-e_3-e_4$. We also note that the lower bound for $m$ is always less
than or equal to the upper bound, which is $2p+1-e_3-e_4$. We thus conclude 
that there are $2p+1-e_3-e_4$ admissible covers with bad reduction.

We now consider the case of an admissible cover with $\rho$ an 
$e_1$-$e_2$-cycle (Theorem
\ref{degenerationlem}.(b)). Let $f_0:V_0 \to X_0$ be such an admissible
cover in characteristic $0$, as
above. In particular, the restriction $f_0^1$ (resp.\ $f_0^2$) has
type $(d_1; e_1, e_2, e_1\text{-}e_2)$ (resp.\ $(d_2; e_1\text{-}e_2, e_3,
e_4)$). 
We write $g_0$ for the Galois closure of $f_0$, and $g_0^i$ for the
corresponding restrictions. Let $G^i$ be the Galois group of $g_0^i$.
The assumptions on the $e_i$ imply that $p$ does not divide the order
of Galois group of $g_0^1$, therefore $g_0^1$ has good reduction to
characteristic $p$.  Moreover, the cover $g_0^1$ is uniquely
determined by the triple $(\rho^{-1}, g_3, g_4)$. If $e_1\neq
e_2$, the gluing is likewise uniquely determined, while if $e_1=e_2$ there 
are exactly $2$ possibilities for the tuple $(g_1, g_2, g_3, g_4)$ for a 
given triple $(\rho^{-1}, g_3, g_4)$.  Therefore to count the number of
admissible covers with bad reduction in this case, it suffices to
consider the reduction behavior of the cover $g_0^2:Y_0^2\to
X_0^2$. 

Corollary \ref{cor:2cyclebad} implies that whether or not $e_1$ equals
$e_2$, the number of admissible covers with bad reduction in the 
$2$-cycle case is equal to $(p+1-e_1-e_2)$ unless $e_1+e_2$ and $e_3$ are 
both even, and bounded from above
by $2(p+1-e_1-e_2)$ always.
  We conclude using Theorem \ref{degenerationlem} that the total
number of admissible covers with bad reduction counted with
multiplicity is less than or equal to
\[
(2p+1-e_3-e_4)+2(p+1-e_1-e_2)=p+(p+1-e_1-e_2)<2p,
\]
and equal to
\[
(2p+1-e_3-e_4)+(p+1-e_1-e_2)=p
\]
unless $e_1+e_2$ and $e_3$ are both even.
The proposition follows.
\end{proof}

\begin{rem} Theorem 4.2 of \cite{os7} does not need the assumption $d=p$. 
Therefore the proof of Proposition \ref{prop:admbad} in the
single-cycle case shows the following stronger result. Let $(d; e_1,
e_2, e_3, e_4)$ be a genus-$0$ type with $1<e_1\leq e_2\leq e_3\leq
e_4 < p$. Then the number of admissible covers with a single ramified
point over the node and bad reduction to characteristic $p$ is
$$(d-p+1)(d+p+1-e_3-e_4)$$
when either $d+1 \geq e_2+e_3$ or $d+1-e_1 < p$. Otherwise, all
admissible covers have bad reduction.
\end{rem}

\section{Proof of the main result}\label{sec:4pt}
In this section, we count the number of mere covers with
ramification type $(p;e_1, e_2, e_3, e_4)$ and bad reduction in the
case that the branch points are generic. Equivalently, we compute the
$p$-Hurwitz number $h_p(p;e_1, e_2, e_3, e_4)$. 

Suppose that $r=4$ and fix a genus-$0$ type $\tau=(p; e_1, e_2, e_3,
e_4)$ with $2\leq e_1\leq e_2\leq e_3\leq e_4<p$.
  We let $g:Y\to X= \PP^1_K$ be a
Galois cover of type $\tau$ defined over a local field $K$ as in
Section \ref{sec:stable}, such that $(X; Q_i)$ is the generic $r$-marked
curve of genus $0$. It is no restriction to suppose that $Q_1=0,
Q_2=1, Q_3=\lambda, Q_4=\infty$, where $\lambda$ is transcendental
over $\QQ_p$.  We suppose that $g$ has bad reduction to characteristic
$p$, and denote by $\bar{g}:\bar{Y}\to \bar{X}$ the stable
reduction. We have seen in Section \ref{sec:stable} that we may associate
with $\bar{g}$ a set of primitive tail covers $( \bar{g}_i)$ and a
deformation datum $(\bar{Z}_0, \omega)$. The primitive tail
covers $\bar{g}_i$ for $i\in\BB=\{1,2,3,4\}$ are  uniquely
determined by the $e_i$ (Lemma \ref{lem:tail1}).

The following proposition shows that the number of covers with bad
reduction is divisible by $p$ in the case that the branch point are
generic.

\begin{prop} \label{prop:baddeg}
Suppose that $(X=\PP^1_K; Q_i)$ is the generic $r=4$-marked curve of
genus zero. Then the number of mere covers of $X$ of ramification type
$(p;e_1, e_2, e_3, e_4)$ with bad reduction is nonzero and divisible by $p$.
\end{prop}

\begin{proof} 
Since the number of Galois covers and the number of mere covers
differ by a prime-to-$p$ factor, it suffices to prove the proposition
for Galois covers.  The existence portion of the proposition is proved
in \cite{bo4}, Proposition 2.4.1, and the divisibility by $p$ in Lemma
3.4.1 of {\it loc.\ cit.}\ (in a more general setting).  We briefly
sketch the proof, which is easier in our case due to the simple
structure of the stable reduction (Lemma \ref{lem:stablered}).  The
idea of the proof is inspired by a result of \cite{we3}, Section~3.

 We begin by observing that away from the wild branch point
 $\xi_i$, the primitive tail cover $\bar{g}_i$ is tamely
 ramified. Therefore we can lift this cover of affine curves to
 characteristic zero.

Let $\cX_{0}=\PP^1_R$ be equipped with $4$ sections $Q_1=0, Q_2=1,
Q_3=\lambda, Q_4=\infty$, where $\lambda\in R$ is transcendental over
$\ZZ_p$. Then (\ref{eq:Kummer}) defines
an $m$-cyclic cover $\cZ_0\to\cX_0$.  We write $Z\to X$ for its generic
fiber. Proposition \ref{prop:dd} implies the existence of a
deformation datum $(\bar{Z}_0, \omega)$. Associated with the deformation
datum is a character $\chi:\ZZ/m\ZZ\to \FF_p^\times$ defined by
$\chi(\beta)=\beta^\ast z/z\pmod{z}$. The differential form $\omega$
corresponds to a $p$-torsion point $P_0\in J(\bar{Z}_0)[p]_\chi$ on
the Jacobian of $\bar{Z}_0$. See for example \cite{se7} 
 (Here we use that the conjugacy classes $C_i$ are
conjugacy classes of prime-to-$p$ elements. This implies that the
differential form $\omega$ is holomorphic.)

 Since $\sum_{i=1}^4 h_{i}=2m$ and the branch points are generic, we
 have that $J(\bar{Z}_0)[p]_\chi\simeq \ZZ/p\ZZ\times {\boldsymbol
   \mu}_p$ (\cite{bo6}, Proposition 2.9) After
 enlarging the discretely valued field $K$, if necessary, we may
 choose a $p$-torsion point $P\in J(\cZ_0\otimes_R K)[p]_\chi$ lifting
 $P_0$. It corresponds to an \'etale $p$-cyclic cover $W\to Z$. The
 cover $\psi:W\to X$ is Galois, with Galois group
 $N:=\ZZ/p\ZZ\rtimes_\chi\ZZ/m\ZZ$. It is easy to see that $\psi$ has
 bad reduction, and that its deformation datum is $(\bar{Z}_0,
 \omega)$.

By using formal patching (\cite{ra3} or \cite{we1}), one now
checks that there exists a map $g_R:\cY\to \cX$ of stable curves over
$\Spec(R)$ whose generic fiber is a $G$-Galois cover of smooth curves,
and whose special fiber defines the given tails covers and the
deformation datum. Over a neighborhood of the original component $g_R$
is the induced cover $\Ind_N^G \cZ_0\to \cX_0$. Over the tails, the
cover $g_R$ is induced by the lift of the tail covers. The fact that
we can patch the tail covers with the cover over $\cX_0$ follows
from the observation that $h_{i}<m_{i}$ (Lemma \ref{lem:tail1}), since
locally there a unique cover with this ramification (\cite{we1}, Lemma
2.12). This proves the existence statement.

The divisibility by $p$ now 
follows from the observation that the set of lifts $P$ of the
$p$-torsion point $P_0\in J(\bar{Z}_0)[p]_\chi$ corresponding to the
deformation datum is a ${\boldsymbol \mu}_p$-torsor.
\end{proof}

We are now ready to prove our Theorem \ref{thm:main}, as well as
a slightly sharper version of Theorem \ref{thm:good-degen}.

\begin{thm}\label{thm:main2} Let  $p$ be an odd prime and $k$ an algebraically
closed field of characteristic $p$. Suppose we are given
integers $2\leq e_1\leq e_2\leq e_3\leq e_4<p$. There exists a dense
open subset $U\subset \PP^1_k$ such that for $\lambda\in U$ the
number of degree-$p$ covers with ramification type $(e_1, e_2, e_3, e_4)$
over the branch points $(0,1,\lambda,\infty)$ is given by the formula
$$h_p(e_1,\dots,e_4)=\min_i(e_i(p+1-e_i))-p.$$ 
Furthermore, unless both $e_1+e_2$ and $e_3$ are even, every such
cover has good degeneration under a degeneration of the base sending
$\lambda$ to $\infty$.
\end{thm}

\begin{proof} Proposition \ref{prop:baddeg} implies that the 
number of covers with ramification type $(p;e_1, e_2, e_3, e_4)$ and
bad reduction is at least $p$. This implies that the generic Hurwitz
number $h_p(e_1,\ldots, e_4)$ is at most $\min_i(e_i(p+1-e_i))-p$.
Proposition \ref{prop:admbad} implies that the number of admissible
covers in characteristic $p$ strictly larger than
$\min_i(e_i(p+1-e_i))-2p$. Since the number of separable covers can
only decrease under specialization, we conclude that the generic
Hurwitz number  equals $\min_i(e_i(p+1-e_i))-p$.  This proves the
first statement, and the second follows immediately from Proposition
\ref{prop:admbad} in the situation that $e_1+e_2$ and $e_3$ are not
both even.
\end{proof}

\begin{rem} By using the results of \cite{bo4} one can  prove a
 stronger result than Theorem \ref{thm:main2}. We say that a
 $\lambda\in \PP^1_k\setminus\{0,1,\infty\}$ is {\em supersingular} if
   it is a zero of the polynomial (\ref{eq:Hasseinv}) and {\em
     ordinary} otherwise. Then the number of covers in characteristic
   $p$ of type $(p; e_1, e_2, e_3, e_4)$ branched at $(0,1,\lambda,
   \infty)$ is $h_p(p; e_1, e_2, e_3, e_4)$ if $\lambda$ is ordinary
   and $h_p(p; e_1, e_2, e_3, e_4)-1$ if $\lambda$ is
   supersingular. To prove this result, one needs to study the stable
   reduction of the cover $\pi:\bar{\cH}\to \PP^1_\lambda$ of the
   Hurwitz curve to the configuration space. We do not prove this
   result here, as it would require too many technical details.
\end{rem}

\bibliographystyle{hamsplain}
\bibliography{hgen}

\end{document}